\documentclass[conference]{IEEEtran}
\usepackage{amsmath}
\usepackage{amssymb}
\usepackage{graphicx}

\ifCLASSINFOpdf
\else
\fi
\begin{document}
%
% paper title
% Titles are generally capitalized except for words such as a, an, and, as,
% at, but, by, for, in, nor, of, on, or, the, to and up, which are usually
% not capitalized unless they are the first or last word of the title.
% Linebreaks \\ can be used within to get better formatting as desired.
% Do not put math or special symbols in the title.
\title{Non-Negative Matrix Factorization Test Cases}

% author names and affiliations
% use a multiple column layout for up to three different
% affiliations
\author{\IEEEauthorblockN{
  Connor Sell and Jeremy Kepner \\
  Massachusetts Institute of Technology, Cambridge, MA 02139 \\
  email: csell@mit.edu, kepner@ll.mit.edu
}
}

%\author{\IEEEauthorblockN{Connor Sell}
%\IEEEauthorblockA{Lincoln Laboratory\\
%Massachusetts Institute of Technology\\
%Cambridge, Massachusetts 02139\\
%Email: csell@mit.edu}
%\and
%\IEEEauthorblockN{Dr. Jeremy Kepner}
%\IEEEauthorblockA{Lincoln Laboratory\\
%Cambridge, Massachusetts 02139\\
%Email: kepner@ll.mit.edu}}

% conference papers do not typically use \thanks and this command
% is locked out in conference mode. If really needed, such as for
% the acknowledgment of grants, issue a \IEEEoverridecommandlockouts
% after \documentclass

% for over three affiliations, or if they all won't fit within the width
% of the page, use this alternative format:
% 
%\author{\IEEEauthorblockN{Michael Shell\IEEEauthorrefmark{1},
%Homer Simpson\IEEEauthorrefmark{2},
%James Kirk\IEEEauthorrefmark{3}, 
%Montgomery Scott\IEEEauthorrefmark{3} and
%Eldon Tyrell\IEEEauthorrefmark{4}}
%\IEEEauthorblockA{\IEEEauthorrefmark{1}School of Electrical and Computer Engineering\\
%Georgia Institute of Technology,
%Atlanta, Georgia 30332--0250\\ Email: see http://www.michaelshell.org/contact.html}
%\IEEEauthorblockA{\IEEEauthorrefmark{2}Twentieth Century Fox, Springfield, USA\\
%Email: homer@thesimpsons.com}
%\IEEEauthorblockA{\IEEEauthorrefmark{3}Starfleet Academy, San Francisco, California 96678-2391\\
%Telephone: (800) 555--1212, Fax: (888) 555--1212}
%\IEEEauthorblockA{\IEEEauthorrefmark{4}Tyrell Inc., 123 Replicant Street, Los Angeles, California 90210--4321}}

% use for special paper notices
%\IEEEspecialpapernotice{(Invited Paper)}

% make the title area
\maketitle

% As a general rule, do not put math, special symbols or citations
% in the abstract
\begin{abstract}
Non-negative matrix factorization (NMF) is a problem with many applications, ranging from facial recognition to document clustering.  However, due to the variety of algorithms that solve NMF, the randomness involved in these algorithms, and the somewhat subjective nature of the problem, there is no clear ``correct answer'' to any particular NMF problem, and as a result, it can be hard to test new algorithms.  This paper suggests some test cases for NMF algorithms derived from matrices with enumerable exact non-negative factorizations and perturbations of these matrices.  Three algorithms using widely divergent approaches to NMF all give similar solutions over these test cases, suggesting that these test cases could be used as test cases for implementations of these existing NMF algorithms as well as potentially new NMF algorithms.  This paper also describes how the proposed test cases could be used in practice.
\end{abstract}

% no keywords

% For peer review papers, you can put extra information on the cover
% page as needed:
% \ifCLASSOPTIONpeerreview
% \begin{center} \bfseries EDICS Category: 3-BBND \end{center}
% \fi
%
% For peerreview papers, this IEEEtran command inserts a page break and
% creates the second title. It will be ignored for other modes.
\IEEEpeerreviewmaketitle

\section{Introduction}
% no \IEEEPARstart
\let\thefootnote\relax\footnotetext{This material is based in part upon work supported by the NSF under grant number DMS-1312831.  Any opinions, findings, and conclusions or recommendations expressed in this material are those of the authors and do not necessarily reflect the views of the National Science Foundation.} What do document clustering, recommender systems, and audio signal processing have in common?  All of them are problems that involve finding patterns buried in noisy data.  As a result, these three problems are common applications of algorithms that solve non-negative matrix factorization, or NMF \cite{Gemulla:2011:LMF:2020408.2020426,NIPS2005_2757,wang2010instantaneous}.

Non-negative matrix factorization involves factoring some matrix $ \mathbf{A} $, usually large and sparse, into two factors $ \mathbf{W} $ and $ \mathbf{H} $, usually of low rank
\begin{equation}
\mathbf{A} = \mathbf{WH}
\end{equation}
Because all of the entries in $ \mathbf{A} $, $ \mathbf{W} $, and $ \mathbf{H} $ must be non-negative, and because of the imposition of low rank on $ \mathbf{W} $ and $ \mathbf{H} $, an exact factorization rarely exists.  Thus NMF algorithms often seek an approximate factorization, where $ \mathbf{WH} $ is close to $ \mathbf{A} $.  Despite the imprecision, however, the low rank of $ \mathbf{W} $ and $ \mathbf{H} $ forces the solution to describe $ \mathbf{A} $ using fewer parameters, which tends to find underlying patterns in $ \mathbf{A} $.   These underlying patterns are what make NMF of interest to a wide range of applications.

In the decades since NMF was introduced by Seung and Lee \cite{NIPS2000_1861}, a variety of algorithms have been published that compute NMF \cite{ALS}.  However, the non-deterministic nature of these NMF algorithms make them difficult to test.  First, NMF asks for approximations rather than exact solutions, so whether or not an output is correct is somewhat subjective.  Although cost functions can quantitatively indicate how close a given solution is to being optimal, most algorithms do not claim to find the globally optimal solution, so whether or not an algorithm gives useful solutions can be ambiguous.  Secondly, all of the algorithms produced so far are stochastic algorithms, so running the algorithm on the same input multiple times can give different outputs if they use different random number sequences.  Thirdly, the algorithms themselves, though often simple to implement, can have very complex behavior that is difficult to understand.  As a result, it can be hard to determine whether a proposed algorithm really ``solves'' NMF.

This paper proposes some test cases that NMF algorithms should solve verifiably.  The approach uses very simple input, such as matrices that have exact non-negative factorizations, that reduce the space of possible solutions and ensure that the algorithm finds correct patterns with little noise.  In addition, small perturbations of these simple matrices are also used, to ensure that small variations in the matrix $ \mathbf{A} $ do not drastically change the generated solution.

% You must have at least 2 lines in the paragraph with the drop letter
% (should never be an issue)

\section{Perturbations of order $ \epsilon $}
Suppose NMF is applied to a non-negative matrix $ \mathbf{A} $ to get non-negative matrices $ \mathbf{W} $ and $ \mathbf{H} $ such that $ \mathbf{A} \approx \mathbf{W} \mathbf{H} $.  If $ \mathbf{A} $ is chosen to have an exact non-negative factorization, then the optimal solution satisfies $ \mathbf{A} = \mathbf{WH} $. Furthermore, if $ \mathbf{A} $ is simple enough, most ``good'' NMF algorithms will find the exact solution.

For example, suppose $ \mathbf{A}_0 $ is a non-negative square diagonal matrix, and the output $ \mathbf{W}_0 $ and $ \mathbf{H}_0 $ is also specified to be square.   Let the diagonal $ n \times n $ matrix $ \mathbf{A}_0 $ be denoted $ \mathbf{A}_0 = \text{diag}(\mathbf{a}_0) $, where $ \mathbf{a}_0 $ is an $ n $-dimensional vector, so that the diagonal entries $ \mathbf{A}_0(i,i) $ are $ \mathbf{a}_0(i) $.  It is easy to show that $ \mathbf{W}_0 $ and $ \mathbf{H}_0 $ must be monomial matrices (diagonal matrices under a permutation) \cite{proof}.  Ignoring the permutation and similarly denoting $ \mathbf{W}_0 = \text{diag}(\mathbf{w}_0) $ and $ \mathbf{H}_0 = \text{diag}(\mathbf{h}_0) $, then $ \mathbf{a}_0(i) = \mathbf{w}_0(i) \mathbf{h}_0(i) $ for applicable $ i $.  Such diagonal matrices $ \mathbf{A}_0 $ were given as input to the known NMF algorithms described in the next section, and all of the algorithms successfully found exact solutions in the form of monomial matrices for $ \mathbf{W}_0 $ and $ \mathbf{H}_0 $.

One way to analyze the properties of an algorithm is to perturb the input by a small amount $ \mathbf{\epsilon} > 0 $ and see how the output changes.  Formally, if the input $ \mathbf{A}_0 $ gives output $ \mathbf{W}_0\mathbf{H}_0 $, then the output generated from $ \mathbf{A}_0 + \epsilon \mathbf{A}_1 $ can be approximated as $ (\mathbf{W}_0 + \epsilon \mathbf{W}_1)(\mathbf{H}_0 + \epsilon \mathbf{H}_1) $.  It is assumed that $ \epsilon $ is sufficiently small that $ \epsilon^2 $ terms are negligible.

For the test case, the nonzero entries of $ \mathbf{A}_1 $ were chosen to be the on the superdiagonal (the first diagonal directly above the main diagonal).  This matrix is denoted  as $ \mathbf{A}_1 = \text{diag}(\mathbf{a}_1, 1) $, where $ \mathbf{a}_1 $ is an $ n-1 $-dimensional vector such that $ \mathbf{A}_1(i,i+1) = \mathbf{a}_1(i) $.  The resulting matrix $ \mathbf{A}_0 + \epsilon \mathbf{A}_1 $ has $O(1)$ entries on its main diagonal, $O(\epsilon) $ entries on the superdiagonal, and zeroes elsewhere.  It is assumed that all the vector entries $ \mathbf{a}_0(i) $ and $ \mathbf{a}_1(i) $ are of comparable magnitude.

\section{Results from Various Algorithms}
Three published NMF algorithms were implemented and run with input of the form $ \mathbf{A} = \mathbf{A}_0 + \epsilon \mathbf{A}_1 $ as described above.  Algorithm 1 was the multiplicative update algorithm described by Seung and Lee in their groundbreaking paper \cite{NIPS2000_1861}, which was run for $ 10^6 $ iterations in each test.  Algorithm 2 was the ALS algorithm described in \cite{ALS}, and which was run for $ 10^6 $ iterations as well.  Algorithm 3 was a gradient descent method as described by Guan and Tao \cite{GradDes}, and was run for $ 10^4 $ iterations.  These three algorithms were chosen because they were representative and easy-to-implement algorithms of three distinct types.  Many published NMF algorithms are variations of these three algorithms.

The experiments began with the simplest nontrivial case, in which $ \mathbf{A} $ is a $ 2 \times 2 $ matrix with only three nonzero entries, with fixed $ \mathbf{a}_0 = [1~1] $ and $ \mathbf{a}_1 = [1] $, while  $ \epsilon $ was varied over several different values.  Each of the algorithms used randomness in the form of initial seed values for $ \mathbf{W} $ and $ \mathbf{H} $. The random seeds were held constant as $ \epsilon $ varied.  As a result, the outputs from the algorithms with different values of $ \epsilon $ were comparable within each test case.

For the $ 2 \times 2 $ case, it is possible to enumerate all of the non-negative exact factorizations of $ \mathbf{A} $.  Given that the factors $ \mathbf{W} $ and $ \mathbf{H} $ are also $ 2 \times 2 $ matrices, they can be written as shown below.
\begin{equation}
\left[ \begin{array}{cc} m & n \\ p & q \end{array} \right] \left[ \begin{array}{cc} r & s \\ t & u \end{array} \right] = \left[ \begin{array}{cc} 1 & \epsilon \\ & 1 \end{array} \right]
\end{equation}
Multiplying the matrices directly produces the the following four equations:
\begin{align}
mr + nt & = 1 \\
ms + nu & = \epsilon \\
pr + qt & = 0 \\
ps + qu & = 1
\end{align}
Recall that all entries must be non-negative, so from equation (5), either $ p $ or $ r $ must be 0, and either $ q $ or $ t $ must be 0.  Furthermore, it cannot be that $ p = q = 0 $ because that would contradict equation (6), and it cannot be that $ r = t = 0 $ because that would contradict equation (3).  Thus two cases remain: $ p = t = 0 $ and $ q = r = 0 $.

Substituting $ p = t = 0 $ into equations (3), (4), and (6) and solving for $ r $, $ s $, and $ u $ gives
\begin{equation}
r = \frac{1}{m}, ~~~
s = \frac{1}{m}\left( \epsilon - \frac{n}{q} \right), ~~~
u = \frac{1}{q}
\end{equation}
Likewise, substituting $ q = r = 0 $ into (3), (4), and (6) and solving for $ s $, $ t $, and $ u $ to gives
\begin{align}
s = \frac{1}{p}, ~~~
t = \frac{1}{n}, ~~~
u = \frac{1}{n} \left( \epsilon - \frac{m}{p} \right)
\end{align}
Observe that these two solutions look similar.  In fact, they differ merely by a permutation.  In the first case,  $ \mathbf{W} $ and $ \mathbf{H} $ have the same main diagonal and superdiagonal format as $ \mathbf{A} $, and can be written in matrix notation as
\begin{equation}
\mathbf{WH} = \left[ \begin{array}{cc} \mathbf{w}_0(1) & \mathbf{w}_1(1) \\ & \mathbf{w}_0(2) \end{array} \right] \left[ \begin{array}{cc} \frac{1}{\mathbf{w}_0(1)} & \frac{1}{\mathbf{w}_0(1)}(\epsilon - \frac{\mathbf{w}_1(1)}{\mathbf{w}_0(2)}) \\ & \frac{1}{\mathbf{w}_0(2)} \end{array} \right]
\end{equation}
The second case can be written as $ (\mathbf{WP})(\mathbf{P}^{-1}\mathbf{H}) $, where $ \mathbf{P} $ is the permutation matrix $ \left[ \begin{array}{cc} & 1 \\ 1 & \end{array} \right] $.

All three of the algorithms tested gave solutions of this form 1000 times out of 1000, for each of several values of $ \epsilon $.  The consistency of the solutions enabled further analysis.  The change in the solution can be measured by the change in the three parameters $ \mathbf{w}_0(1) $, $ \mathbf{w}_0(2) $, and $ \mathbf{w}_1(1) $ (ignoring the permutation if present).  Figure 1 shows the change in each of the three parameters from the base case $ \mathbf{A}_0 $ for several different values of $ \epsilon $ when input into Algorithm 1.  Each of the values is the arithmetic mean of the corresponding values generated from 1000 different random seeds.
\begin{figure}
\begin{center}
\includegraphics[scale=0.22]{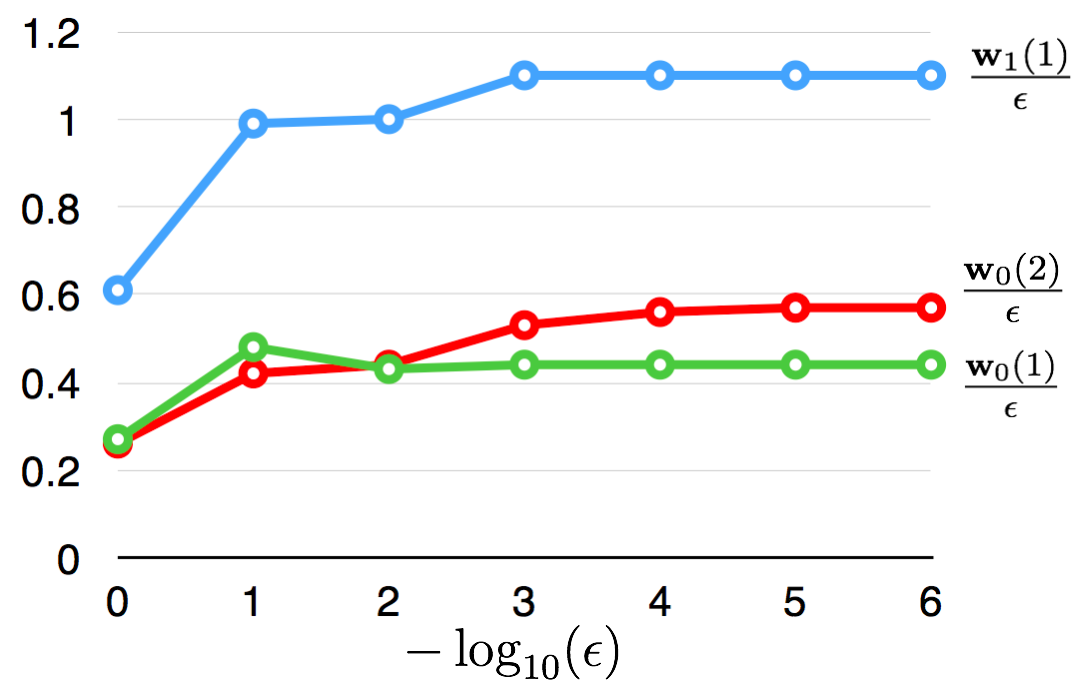}
\caption{The figure shows the slope associated with the change in each of the three parameters for each of several values of $ \epsilon $.  As $ \epsilon $ approaches zero on the right of the graph, the values of the slopes converge, showing that for sufficiently small $ \epsilon $, each of the parameters is linear in $ \epsilon $.}
\end{center}
\end{figure}
Of course, the precise values depend on the distribution of randomness used.  But notice that as $ \epsilon $ approaches 0, the values of the three parameters become very nearly linear in $ \epsilon $.  The results for Algorithms 2 and 3 were very similar - they also showed linearity of the parameters in $ \epsilon $, with comparable slopes.

However, $ \mathbf{w}_1(1) $ was not always linear in $ \epsilon $, even for small $ \epsilon $.  In some cases, the difference approached 0 much more quickly.  To see why this occurred, consider that the entries in $ \mathbf{H} $ could have been chosen to be the parameters rather than the entries in $ \mathbf{W} $.  Also, recall that in the base case $ \mathbf{A}_0 $, in which $ \epsilon = 0 $, $ \mathbf{w}_1(1) = \mathbf{h}_1(1) = 0 $ since both entries are off the diagonal.  Thus, when either is linear in $ \epsilon $, they are of the form $ \epsilon x $ for some slope $ x $.  Since the solution is exact, it can be deduced that
\begin{equation}
\mathbf{w}_0(1) \mathbf{h}_1(1) + \mathbf{w}_1(1) \mathbf{h}_0(2) = \epsilon
\end{equation}
Therefore, in the cases that $ \mathbf{w}_1(1) $ approaches 0 very quickly, since $ \mathbf{w}_0(1) $ approaches a large, stable value as $ \epsilon $ approaches 0, $ \mathbf{h}_1(1) $ must be nearly linear in $ \epsilon $.  So in the cases that $ \mathbf{w}_1(1) $ is not linear in $ \epsilon $, its symmetrical counterpart, $ \mathbf{h}_1(1) $, is.  To simplify this complication out of the data, the parameters in $ \mathbf{W} $ were chosen when $ \mathbf{w}_1(1) $ was closer to linearity in $ \epsilon $, and the parameters in $ \mathbf{H} $ were chosen when $ \mathbf{h}_1(1) $ was closer to linearity in $ \epsilon $.

Curiously, although it was possible for $ \mathbf{w}_1(1) $ and $ \mathbf{h}_1(1) $ to ``split'' the nonlinearity so that both were somewhat linear, this rarely occurred.  All three algorithms preferred to make one of them very close to linear at the expense of the other.  When $ \mathbf{w}_1(1) $ approached zero very rapidly, by equations (3) and (4), $ \mathbf{h}_1(1) = \epsilon \mathbf{h}_0(1) $, and similarly, when $ \mathbf{h}_1(1) $ is negligible, $ \mathbf{w}_1(1) = \epsilon \mathbf{h}_0(2) $.

Next, different values for the entries of $ \mathbf{a}_0 $ and $ \mathbf{a}_1 $ were tried, so they had a range of entries rather than all 1's. The algorithms all behaved similarly; up to permutation, they satisfied the following formula
\begin{equation}
\mathbf{WH} = \left[ \begin{array}{cc} \mathbf{w}_0(1) & \mathbf{w}_1(1) \\ & \mathbf{w}_0(2) \end{array} \right] \left[ \begin{array}{cc} \frac{\mathbf{a}_0(1)}{\mathbf{w}_0(1)} & \frac{\mathbf{a}_1(1)}{\mathbf{w}_0(1)}(\epsilon - \frac{\mathbf{w}_1(1)\mathbf{a}_0(2)}{\mathbf{a}_1(1)\mathbf{w}_0(2)}) \\ & \frac{\mathbf{a}_0(2)}{\mathbf{w}_0(2)} \end{array} \right]
\end{equation}
Note that equation (9) is just a special case of this equation in which $ \mathbf{a}_0(1) = \mathbf{a}_0(2) = \mathbf{a}_1(1) = 1 $.   The same phenomena was also observed in which the algorithm usually made one of $ \mathbf{w}_1(1) $ and $ \mathbf{h}_1(1) $ be nearly linear in $ \epsilon $ and the other approach zero rapidly, rather than having both entries be non-negligible.  As long as the entries of $ \mathbf{a}_0 $ and $ \mathbf{a}_1 $ are roughly on the order of 1, the algorithms operated similarly.

The next case examined set $ \mathbf{A} $ to be a $ 3 \times 3 $ matrix.  Using similar logic to the $ 2 \times 2 $ case, it can be deduced that any exact factorization of $ \mathbf{A} $ is likely to be of the form
\begin{equation}
\left[ \begin{array}{ccc} \mathbf{w}_0(1) & \mathbf{w}_1(1) & \\ & \mathbf{w}_0(2) & \mathbf{w}_1(2) \\ & & \mathbf{w}_0(3) \end{array} \right] \left[ \begin{array}{ccc} \mathbf{h}_0(1) & \mathbf{h}_1(1) & \\ & \mathbf{h}_0(2) & \mathbf{h}_1(2) \\ & & \mathbf{h}_0(3) \end{array} \right]
\end{equation}
Indeed, all three algorithms always gave solutions of this form.  In fact, most of the time there were two more zero entries than necessary - either $ \mathbf{w}_1(1) $ or $ \mathbf{h}_1(1) $, and either $ \mathbf{w}_1(2) $ or $ \mathbf{h}_1(2) $.  This is similar to the way that $ \mathbf{w}_1(1) $ or $ \mathbf{h}_1(1) $ often approached 0 rapidly in the $ 2 \times 2 $ case.  To note another similarity to the $ 2 \times 2 $ case, whenever $ \mathbf{w}_1(i) $ was significant and $ \mathbf{h}_1(i) $ was not, $ \mathbf{w}_1(i) $ was very close to $ \epsilon \mathbf{w}_0(i+1) $ - in similar situations $ \mathbf{h}_1(i) $ was approximately $ \epsilon \mathbf{h}_0(i) $.

As a result, there were 4 distinct configurations of the nonzero elements in the solutions, as given by Figure 2.
\begin{figure}
\begin{center}
\begin{tabular}{|c|c|} \hline
Type & equal to 0 \\ \hline
Type I & $ \mathbf{w}_1(1), \mathbf{w}_1(2) $ \\ \hline
Type II & $ \mathbf{h}_1(1), \mathbf{h}_1(2) $ \\ \hline
Type III & $ \mathbf{w}_1(1), \mathbf{h}_1(2) $ \\ \hline
Type IV & $ \mathbf{h}_1(1), \mathbf{w}_1(2) $ \\ \hline
\end{tabular}
\caption{We categorized the solutions when $ \mathbf{A} $ was a $ 3 \times 3 $ matrix by where the non-negligible entries in the solution were.  For each type, this table shows which entries that are usually positive are negligible.}
\end{center}
\end{figure}
Note that Type IV appears to be an inexact solution; since it has positive $ \mathbf{w}_1(1) $ and $ \mathbf{h}_1(2) $, the entry at position $ \mathbf{A}(1,3) = \mathbf{w}_1(1)\mathbf{h}_1(2) $ in the product $ WH $ would have to be nonzero.  However, both $ \mathbf{w}_1(1) $ and $ \mathbf{h}_1(2) $, like all entries on the superdiagonal, are $ O(\epsilon) $, so $ \mathbf{w}_1(1) \mathbf{h}_1(2) $ is $ O(\epsilon^2) $, and is considered  negligible.  In fact, most of the solutions generated by the algorithms had nonzero values for entries that were supposed to be zero, but for this analysis anything below $ O(\epsilon^2) $ was considered negligible.
%We suspect that these negligible values arise because we cut off the algorithms before they were fully able to converge, as that would take a very long time.

Each algorithm was run 100 times on the $ 3 \times 3 $ input with $ \mathbf{w}_0 = [1~1~1] $, $ \mathbf{w}_1 = [1~1] $, and $ \epsilon = 10^{-3} $.  The solutions were categorized by the solution type in Figure 2.  The distributions of the solutions by algorithm type are given in Figure 2.  Note that some solutions did not have two negligible entries among $ \mathbf{w}_1(1) $, $ \mathbf{w}_1(2) $, $ \mathbf{h}_1(1) $, and $ \mathbf{h}_1(2) $, in which case the smaller entry was ignored for the sake of sorting - this accounted for about 20\% of the three algorithms, the majority occurring in Algorithm 1.
\begin{figure}
\begin{center}
\includegraphics[scale=0.45]{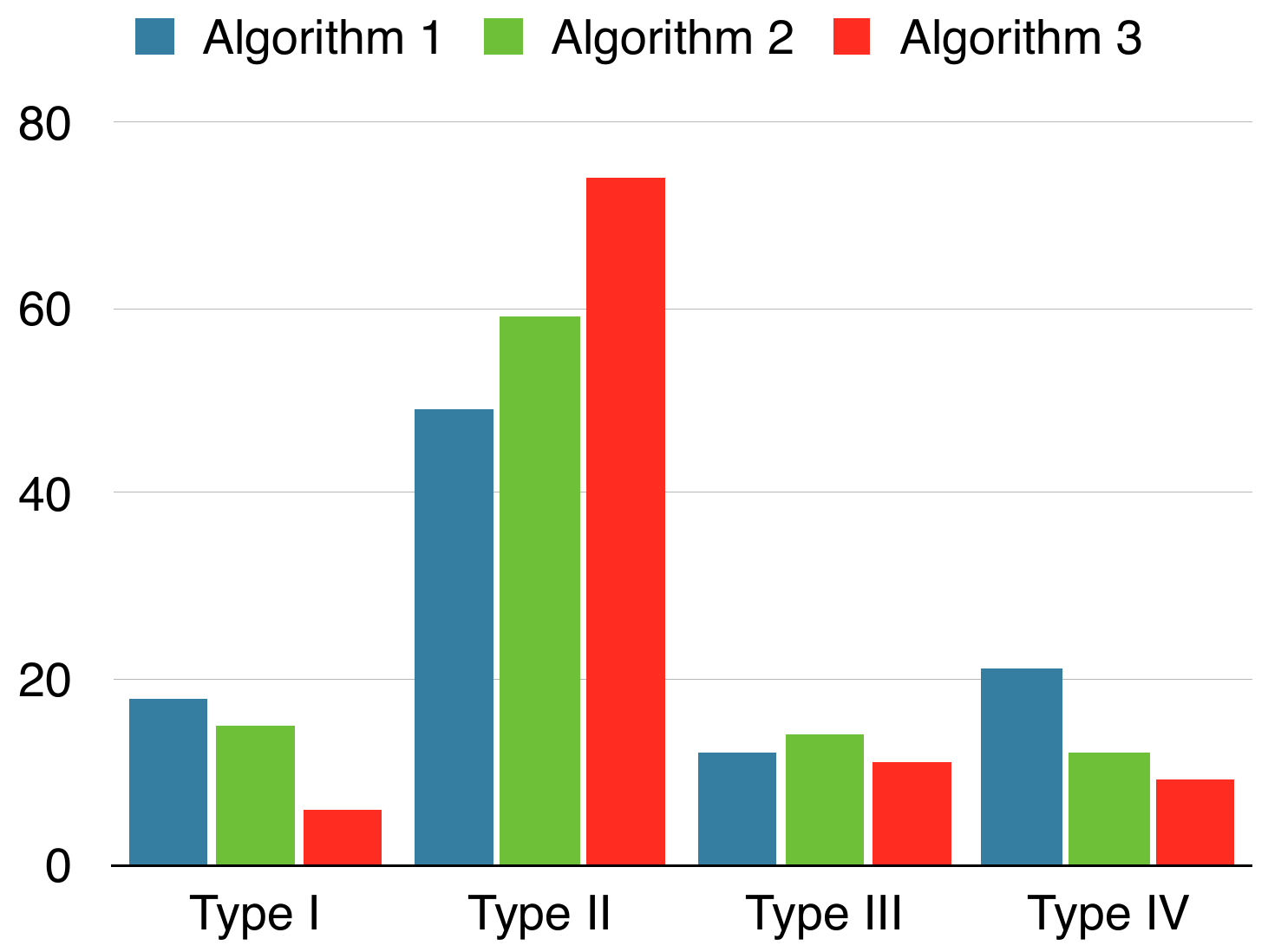}
\caption{Categorized the solutions for $ \mathbf{A} $  being a $ 3 \times 3 $ matrix by where the non-negligible entries in the solution were.  This chart shows how often each algorithm generated a solution of each type out of 100 cases.  Type II (in which $ \mathbf{H} $ is diagonal) was the most common among all the algorithms, but by differing amounts.}
\end{center}
\end{figure}
It is significant to note that even the solutions that didn't fall cleanly into a ``type'' still satisfied the pattern given in (12).  It seems that an NMF algorithm should satisfy this pattern, but little more is required.

Next, entries in $ \mathbf{a}_0 $ and $ \mathbf{a}_1 $, were changed as in the $ 2 \times 2 $ case.  As long as the entries were $O(1)$  (as opposed to $ O(\epsilon) $ or $ O(\frac{1}{\epsilon}) $), the behavior of the algorithms was similar.

Finally, $ \mathbf{A} $ larger than $ 3 \times 3 $ were examined.  Several different sizes of matrices were tested, ranging from $ 4 \times 4 $ to $ 20 \times 20 $, always keeping $ \mathbf{A} $, $ \mathbf{W} $, and $ \mathbf{H} $ square, with positive entries only on the main diagonal and the superdiagonal.  The experiments followed the same general pattern; nonzero entries in $ \mathbf{W} $ and $ \mathbf{H} $ appeared only on the main diagonal and superdiagonal.  Using similar logic to the $ 2 \times 2 $ and $ 3 \times 3 $ cases, it can be shown that these are the only exact solutions.  However, in practice, as the matrices get larger, exceptions to this pattern become more common, particularly in Algorithm 3.  The general rule seems to mostly hold (over half the time) until $ \mathbf{A} $ becomes around $ 20 \times 20 $.  Note, however, that because the run-time of the algorithms are cubic in the size of the matrix, at best, the sample size for large matrices is small.

\section{Proposed Tests for NMF Algorithms}
Since all three algorithms, which cover a variety of approaches to NMF, had a lot in common in their solutions, it is propose that these inputs $ \mathbf{A} $ could be used as a test case of an NMF algorithm implementation.  In this section, it is proposed how such test cases could be executed.

The test begins with input of the form
\begin{equation}
\mathbf{A} = \mathbf{A}_0 + \epsilon \mathbf{A}_1 = \text{diag}(\mathbf{a}_0) + \epsilon \text{diag}(\mathbf{a}_1,1)
\end{equation}
$ \mathbf{A} $ is square, and preferably somewhere between $ 3 \times 3 $ and $ 8 \times 8 $ in size, although bigger inputs may be useful as well.  The entries should vary between tests.  Each test should start by using $ \epsilon = 0 $ so that $ \mathbf{A} $ is diagonal.  The results of this test should have $ \mathbf{W} $ and $ \mathbf{H} $ monomial - only one nonzero element in each row and column.  Ignore entries that are below $ O(10^{-10}) $, for the entirety of testing, as any such entries are negligible.

If $ \mathbf{W} $ or $ \mathbf{H} $ is not monomial, or if the product $ \mathbf{WH} $ is not equal to $ \mathbf{A} $ to within a negligible margin of error, the algorithm fails the test.  Otherwise, the generated solution can be used to find the permutation matrix $ \mathbf{P} $ that makes $ \mathbf{WP} $ and $ \mathbf{P}^{-1}\mathbf{H} $ diagonal by replacing the nonzero entries of $ \mathbf{H} $ with 1's.  Since $ \mathbf{A} = \mathbf{WH} $ is diagonal, $ \mathbf{WP} $ is also diagonal, and since $ \mathbf{I} = \mathbf{P}^{-1}\mathbf{P} $ is diagonal, so is $ \mathbf{P}^{-1}\mathbf{H} $.  Knowing $ \mathbf{P} $ will make the rest of the testing much simpler since it is easier to identify whether a solution is of the form given above when it is not permuted.

Next, run the test again using a positive value for $ \epsilon $; $ \epsilon = 10^{-3} $ seems to work well, although using a variety of $ \epsilon $ is also recommended.  Make sure to use the same random seeds that were used in the $ \epsilon = 0 $ test to produce corresponding output.  Then check that the $ \mathbf{W} $ and $ \mathbf{H} $ given by the algorithm are such that $ \mathbf{WP} $ and $ \mathbf{P}^{-1}\mathbf{H} $ have nonzero entries only on the two diagonals that they are supposed to.  If this doesn't hold, changing $ \epsilon $ might have changed which permutation returns $ \mathbf{W} $ and $ \mathbf{H} $ to the proper form, so check again; this happens more commonly among larger matrices than smaller ones.  However, if $ \mathbf{W} $ and $ \mathbf{H} $ really do break the form, or $ \mathbf{A} \neq \mathbf{WH} $, the algorithm fails the test on this input.  Otherwise, it passes.

Note that even widely accepted algorithms do fail these tests occasionally, especially with matrices larger than $ 8 \times 8 $, so it's advisable to perform the test many times to get a more accurate idea of an algorithm's performance.

\section{Conclusion}
This paper proposes an approach to the problem of testing NMF algorithms by running the algorithms on simple input that can produce an exact non-negative factorization, and perturbations of such input.  In particular,  square matrices with $O(1)$ entries on the main diagonal and $ O(\epsilon)$ entries on the superdiagonal are proposed, because they have exact solutions that can enumerated mathematically, or because they are perturbations of matrices with exact solutions.
%, and because they're simple enough for useful algorithms to find exact solutions.

The test cases have been used as input on three known NMF algorithms that represent a variety of algorithms, and all of them behaved similarly, which suggests testable, quantifiable behaviors that many NMF algorithms share.  These test cases offer one approach for testing candidate NMF implentations to help determine whether it behaves as it should.

% conference papers do not normally have an appendix

% use section* for acknowledgment
\section*{Acknowledgment}
The authors would like to thank Dr. Alan Edelman for providing and overseeing this research opportunity, and Dr. Vijay Gadepally for his advice and expertise.

% trigger a \newpage just before the given reference
% number - used to balance the columns on the last page
% adjust value as needed - may need to be readjusted if
% the document is modified later
%\IEEEtriggeratref{8}
% The "triggered" command can be changed if desired:
%\IEEEtriggercmd{\enlargethispage{-5in}}

% references section

% can use a bibliography generated by BibTeX as a .bbl file
% BibTeX documentation can be easily obtained at:
% http://mirror.ctan.org/biblio/bibtex/contrib/doc/
% The IEEEtran BibTeX style support page is at:
% http://www.michaelshell.org/tex/ieeetran/bibtex/
%\bibliographystyle{IEEEtran}
% argument is your BibTeX string definitions and bibliography database(s)
%\bibliography{IEEEabrv,../bib/paper}
%
% <OR> manually copy in the resultant .bbl file
% set second argument of \begin to the number of references
% (used to reserve space for the reference number labels box)

% that's all folks
\end{document}